\documentclass{conm-p-l}
\usepackage{amsmath,amssymb}

\newcommand{\tbf}{\textbf}
\newcommand{\tit}{\textit}

\newcommand{\mcal}{\mathcal}

\newcommand{\mbb}{\mathbb}

\newcommand{\CC}{\mbb{C}}

\newcommand{\RR}{\mbb{R}}
\newcommand{\ZZ}{\mbb{Z}}

\newcommand{\calb}{\mcal{B}}
\newcommand{\calc}{\mcal{C}}
\newcommand{\cald}{\mcal{D}}

\newcommand{\calh}{\mcal{H}}
\newcommand{\call}{\mcal{L}}

\newcommand{\calo}{\mcal{O}}
\newcommand{\caln}{\mcal{N}}

\newcommand{\cals}{\mcal{S}}

\newcommand{\al}{\alpha}

\newcommand{\del}{\delta}

\newcommand{\lam}{\lambda}

\newcommand{\uxg}{\U_\chi(\frg)}

\newcommand{\frb}{\mathfrak{b}}

\newcommand{\frg}{\mathfrak{g}}

\newcommand{\frp}{\mathfrak{p}}

\newcommand{\frsl}{\mathfrak{s}\mathfrak{l}}

\DeclareMathOperator{\Lie}{Lie}

\DeclareMathOperator{\U}{U}

\newtheorem*{conj}{Conjecture}

\theoremstyle{definition}

\theoremstyle{remark}

\numberwithin{equation}{section}

\begin{document}

\title[Reduced enveloping algebras and cells in the affine Weyl
group]{Representations of reduced enveloping algebras and cells in the
affine Weyl group}

\author{J.E. Humphreys}
\address{Dept. of Mathematics \& Statistics, U. Massachusetts, Amherst, MA 01003}

\email{jeh@math.umass.edu} \thanks{For helpful advice I am grateful to
Roman Bezrukavnikov, Paul Gunnells, Jens C. Jantzen, Victor Ostrik,
Jian-yi Shi, and Eric Sommers.}

\date{}

\begin{abstract} Let $G$ be a semisimple algebraic group over an
algebraically closed field of characteristic $p>0$, and let
$\mathfrak{g}$ be its Lie algebra.  The crucial Lie algebra
representations to understand are those associated with the reduced
enveloping algebra $U_\chi(\mathfrak{g})$ for a ``nilpotent'' $\chi
\in \mathfrak{g}^*$.  We conjecture that there is a natural assignment
of simple modules in a regular block to left cells in the affine Weyl
group $W_a$ (for the dual root system) lying in the two-sided cell
which corresponds to the orbit of $\chi$ in Lusztig's bijection.  This
should respect the action of the component group of $C_G(\chi)$ and
fit naturally into Lusztig's enriched bijection involving the
characters of $C_G(\chi)$.  Some evidence will be described in special
cases.  \end{abstract}

\maketitle

In order to explain the conjecture, we have to review some facts about
three logically independent topics:
\begin{itemize} 

\item[(A)] cells in affine Weyl groups, 

\item[(B)] nilpotent orbits, 

\item[(C)] Lie algebra representations in characteristic $p>0$.
\end{itemize} Subtle connections between (A) and (B) have been
discovered by Lusztig, while connections between (B) and (C) have
emerged over several decades (notably in the work of Kac--Weisfeiler,
Friedlander--Parshall, Premet, and others cited below).  We hope to
build further links between (A) and (C), with the goal of finding
a representation-theoretic model for Lusztig's formal conjecture in
\cite[\S 10]{lu89}.

Notation varies considerably in the literature (and sometimes
clashes).  Our conventions here start with a simple, simply connected
algebraic group $G$ over an algebraically closed field $K$ of
characteristic $p>0$.  Let $T$ be a maximal torus and $W$ the Weyl
group.  Denote by $\Phi$ the root system, with positive system
$\Phi^+$ relative to a simple system $\Delta$.  The character group
$X=X(T)$ is the full weight lattice for $\Phi$.  Let $Q=\ZZ \Phi$ be
the root lattice.

\section{Cells in affine Weyl groups}

\subsection{}

First we recall some basic results about cells.  These arise in the
work of Kazhdan and Lusztig on arbitrary Coxeter groups and their
Hecke algebras, but here we focus just on the case of affine Weyl
groups.  (See Lusztig's papers \cite{lu83}--\cite{lu89} as well as Shi
\cite{sh86}, Xi \cite{xi94b,xi02}.)

Define $W_a := W \ltimes Q$ (the \emph{affine Weyl group}) and
$\widetilde{W}_a := W \ltimes X$ (the \emph{extended affine Weyl
group}).  The latter is not usually a Coxeter group, but is important
for Lusztig's $p$-adic group program; here we focus just on $W_a$.  We
say $W_a$ is of type $\widetilde{X}_n$ if $\Phi$ is of type $X_n$.  It
is important to note that $W_a$ is a \emph{dual} version of the usual
affine Weyl group constructed by Bourbaki via the coroot lattice; this
reflects the influence of Langlands duality in Lusztig's program.

As in the case of an arbitrary Coxeter group, the group $W_a$ is
partitioned into \emph{two-sided cells} (here denoted $\Omega$).  Each
of these is in turn partitioned into \emph{left cells} (here denoted
$\Gamma$) or equally well into \emph{right cells}, each of which is
the set of inverses of elements in some left cell.  These partitions
arise (together with Kazhdan--Lusztig polynomials) from comparison of
the Kazhdan--Lusztig basis for the Hecke algebra with the standard
basis.

The definition of cells yields a natural partial ordering on the
collection of two-sided cells.  The \emph{highest} cell in this
ordering contains just the identity element $1$ of $W_a$.

Since $W_a$ acts simply transitively on the \emph{alcoves} in the
affine space $E:= \RR \otimes_\ZZ X$, the various cells can be
identified with sets of alcoves.  In this picture $W$ labels the
family of alcoves around the special point $0$.  (Conventions differ
in the literature; for example, some authors work with right actions
rather than left actions in this context.)

\subsection{}

Beyond these generalities, Lusztig develops more special features of
cells for $W_a$.  Generalizing the case of a Weyl group, he defines in
\cite{lu85} an \emph{a-invariant} $a(w)$ for each $w \in W_a$,
constant on each two-sided cell and denoted $a(\Omega)$.  This is an
integer between 0 and $N:=|\Phi^+|$, defined combinatorially in terms
of the Hecke algebra.  The $a$-invariant respects (inversely) the
partial ordering of two-sided cells.  For example, $a(\Omega)=0$
precisely when $\Omega= \{1\}$.  At the other extreme, it turns out
that there is a unique cell $\Omega$ with $a(\Omega)=N$; this is the
lowest two-sided cell \cite{sh88}.

\subsection{}

With the help of the $a$-invariant, Lusztig shows that $W_a$ has only
finitely many two-sided cells, each partitioned into finitely many
left (or right) cells.  It is then natural to ask how many two-sided
and one-sided cells there are.  These questions are extremely
difficult to approach in a purely combinatorial way, though they have
been answered for type $\widetilde{A}_n$ by Shi \cite{sh86} and in
some isolated low rank cases.  To formulate and prove general
conjectures, some connection with the geometry of the nilpotent
variety and flag variety seems to be essential.

\subsection{}

The lowest two-sided cell $\Omega$ has been explored thoroughly by Shi
\cite{sh88}.  It contains $|W|$ left cells, each obtained by
intersecting $\Omega$ with a Weyl chamber.  The entire antidominant
chamber is one left cell.  On the other hand, the intersection
$\Gamma$ of $\Omega$ with the dominant chamber is a shifted version of
this chamber: Taking $0$ as the origin in $E$, consider the special
point $\rho$ (the sum of fundamental weights) at which translates of
all root hyperplanes meet.  Then $\Gamma$ consists of the alcoves
lying on the positive sides of all these hyperplanes.

\subsection{}

In \cite{lu87a}, Lusztig defines a set $\cald$ of \emph{distinguished
involutions} in $W_a$, as follows.  For $w \in W_a$, let $\ell(w)$ be
its length and let $\del(w)$ be the degree of the Kazhdan--Lusztig
polynomial $P_{1,w}(q)$.  Then $w \in \cald$ iff $a(w) =
\ell(w)-2\del(w)$, in which case $w$ is shown to be an involution.
Each left cell contains a unique distinguished involution.  For
example, in the dominant left cell of the lowest two-sided cell
$\Omega$ described above, the distinguished involution belongs to the
lowest of the $|W|$ alcoves around the special point $2\rho$.  The
distinguished involution in the antidominant left cell of $\Omega$ is
the longest element of $W$ (if $1$ corresponds to the lowest dominant
alcove).

\subsection{}

Lusztig and Xi \cite{lx88} show that each two-sided cell of $W_a$
contains a \emph{canonical left cell}, whose corresponding alcoves all
lie in the dominant Weyl chamber $\calc \subset E$.  In this way,
$\calc$ is partitioned into canonical left cells belonging to the
two-sided cells.

Chmutova--Ostrik \cite{co02} develop an algorithm to compute the
distinguished involutions in all canonical left cells, with explicit
tables given in low ranks.  But it seems to be more difficult to
locate these involutions in arbitrary left cells.

\subsection{}

Pictures of the cells for the affine Weyl groups of types
$\widetilde{A}_2, \widetilde{B}_2, \widetilde{G}_2$ are given by
Lusztig \cite[\S11]{lu85}.  Paul Gunnells has used computer graphics
to investigate all three-dimensional cells as well.  

Jian-yi Shi \cite{sh86} has worked out the combinatorics in
considerable detail for type $\widetilde{A}_n$, while developing
general tools such as ``sign types'' for the study of cells.  Other
affine Weyl groups of low rank have been studied in a similar spirit
by him and a number of other people, including Robert B\'edard, Cheng
Dong Chen, Jie Du, Gregory Lawton, Feng Li, Jia Chun Liu, He Bing Rui,
Nanhua Xi, Xin Fa Zhang.

\section{Nilpotent orbits and cells}

\subsection{}

Denote by $\caln$ the set of nilpotent elements in $\frg:=\Lie G$.
This is the \emph{nilpotent variety} (or \emph{nullcone}).  It
consists of finitely many orbits under the adjoint action of $G$,
partially ordered by inclusion of one orbit in the closure of another.
The orbits range from $\{0\}$ to the regular orbit, which is dense in
$\caln$ and therefore has dimension $2N =|\Phi|$.  Whenever $p$ is a
``good'' prime (as in 3.1 below), there is a $G$-equivariant
isomorphism between $\caln$ and the unipotent variety of $G$.
Moreover, the partially ordered set of $G$-orbits in $\caln$ is
isomorphic to the corresponding set for the Lie algebra over $\CC$ of
the same type.  Although Lusztig's use of unipotent classes is based
in characteristic 0, the ideas therefore transfer readily to our
situation.  (Jantzen \cite{ja04b} gives a helpful account with
emphasis on characteristic $p$.)

Various other varieties and groups are associated with $\caln$.  The
flag variety $\calb$ of $G$ may be identified with the collection of
Borel subalgebras of $\frg$.  If $e \in \caln$, the set of Borel
subalgebras containing $e$ is denoted by $\calb_e$.  It plays an
essential role in the Springer resolution of singularities of $\caln$,
where it is referred to as a \emph{Springer fiber}.  

Let $C_G(e)$ be the centralizer of $e$ in $G$, and denote by $A(e)$
the finite \emph{component group} $C_G(e)/C_G(e)^\circ$.  The
cohomology of $\calb_e$ with suitable coefficients (complex or
$l$-adic) vanishes in odd degrees and has commuting actions by the
finite groups $W$ and $A(e)$.  We write simply $H^i(\calb_e)$.  This
is the framework for the Springer Correspondence (see for example
\cite[\S13]{ja04b}).

\subsection{}

Soon after the Kazhdan--Lusztig theory was developed, Lusztig
\cite[3.6]{lu83} conjectured the existence of a bijection between the
collection of two-sided cells of $W_a$ (based as above on the root
lattice rather than coroot lattice) and the collection of unipotent
classes in $G$ (or equivalently, the collection of nilpotent orbits
in $\frg$).  This bijection should respect the natural partial orderings,
with the cell $\{1\}$ corresponding to the regular nilpotent orbit and
the lowest two-sided cell corresponding to the zero orbit.  (His ideas
were formulated in characteristic $0$ but adapt to our setting when
$p$ is good.)

By combining a number of deep techniques, Lusztig was able to
construct a suitable bijection in \cite{lu89}. Under his bijection, if
the two-sided cell $\Omega$ corresponds to the orbit of some $e \in
\caln$, then $a(\Omega) = \dim \calb_e$.  But the order-preserving
property remained elusive except in low ranks.  This was later proved
combinatorially for type $A_n$ by Shi, while the general case follows
from recent work of Bezrukavnikov \cite[Thm.~4]{be02}.

\subsection{}

In \cite[3.6]{lu83}, Lusztig formulated a further conjecture on left
cells in terms of the fixed points of $A(e)$ on the cohomology of
$\calb_e$:
\begin{quote}
(LC) The number of left cells in the two-sided cell corresponding to a
nilpotent $e$ should be equal to $\sum_i (-1)^i \dim
H^i(\calb_e)^{A(e)}$. 
\end{quote}
Due to the vanishing of cohomology in odd degrees, the contributions
here are all nonnegative.

While (LC) has not yet been proved in general, it agrees with direct
calculations in low ranks and with the results of Shi for type
$\widetilde{A}_n$ \cite[14.4.5,15.1,17.4]{sh86}.  Here all component
groups are trivial, while on the other hand the representation of $W$
on the cohomology is known to be induced from the trivial character of
a parabolic subgroup $W_I$ generated by reflections relative to a set
$I$ of simple roots.  (See the discussion in \cite[p.~203]{ja04b}).
When translated into the language of partitions, the number
$|W|/|W_I|$ agrees with the number of left cells found by Shi for a
corresponding two-sided cell.

\subsection{}

As part of his more refined study of the ``asymptotic Hecke algebra''
in connection with $p$-adic representations, Lusztig \cite[\S10]{lu89}
formulated more detailed conjectures relating the cells with geometry.
Fix a two-sided cell $\Omega$ corresponding in his bijection to the
orbit of $e \in \caln$, and let $\Gamma_\Omega$ be its canonical left
cell.  Denote by $F$ a maximal reductive subgroup of $C_G(e)$, so
$F/F^\circ \cong A(e)$.  Write $\widehat{F}$ for the set of
isomorphism classes of irreducible representations of $F$.  

Lusztig's conjectural set-up involves a finite set $Y$, acted on by
$A(e)$, with cardinality equal to the Euler characteristic of
$\calb_e$.  The orbits of $A(e)$ in $Y$ should be in bijection with
the left cells in $\Omega$, with a singleton orbit expected to
correspond to the canonical left cell.  In general, the isotropy group
in $A(e)$ of an element $y \in Y$ corresponds to an intermediate
subgroup $F \supset F_y \supset F^\circ$.  The representations of $F$
or $F_y$ enter via a notion of ``$F$-vector bundle'' on $Y$ or $Y
\times Y$.

This formalism is then subject to several requirements in
\cite[10.5]{lu89}.  For example, the representation of $A(e)$ on
$H^*(\calb_e)$ should be equivalent to the permutation representation
of $A(e)$ on $Y$.  (This recovers the statement (LC) above.)  As a
consequence, one should have a natural bijection between
$\Gamma_\Omega \cap \Gamma_\Omega^{-1}$ and $\widehat{F}$.  For an
arbitrary left cell $\Gamma$ corresponding to the orbit of $y \in Y$,
the group $F$ should be replaced by the group $F_y$.

Out of this abstract framework emerges a conjectural bijection between
pairs $(\calo_e,\varphi)$ and $X^+$, where $\calo_e$ is a nilpotent
orbit and $\varphi$ an irreducible representation of $C_G(e)$.  (Such
a bijection was conjectured independently by Vogan.)  Note that when
we work with $W_a$ rather than $\widetilde{W}_a$, the root lattice $Q$
replaces $X$.  

Bezrukavnikov has found suitable bijections in \cite{be03} and
\cite{be04a}; these are shown in \cite[Remark~6]{be02} to coincide.
For other work related to Lusztig's conjectures (especially in this
last formulation), see the individual and joint papers by Achar and
Sommers \cite{so01,as02,ac04}, Bezrukavnikov and Ostrik
\cite{os00,bo04}, Lusztig \cite{lu97a}, Xi \cite{xi94a,xi94b,xi02}.

\section{Lie algebra representations in characteristic $p>0$}

\subsection{}

The representation theory of $\frg$ has been studied over a long
period of time: for surveys of earlier work, see \cite{go01} and
\cite{hu98}.  In a series of papers, Jantzen \cite{ja98}--\cite{ja04a}
has extended the theory considerably.  Here we focus on just the
\emph{simple} modules for the universal enveloping algebra $U(\frg)$.
These all occur as modules for reduced enveloping algebras $\uxg$,
which are finite-dimensional quotients of $U(\frg)$ parametrized by
$\chi \in \frg^*$.  Those $\uxg$ for $\chi$ in a single orbit under
the coadjoint action of $G$ are isomorphic, so one looks for a
well-chosen orbit representative $\chi$.

In order to obtain uniform results, Jantzen imposes several relatively
weak hypotheses (H1)--(H3) on $\frg$ and $p$, which we also assume.
For a simply connected group, he requires the prime $p$ to be
\emph{good} for $\Phi$, which eliminates some root systems when
$p=2,3,5$.  Moreover, the algebras $\frsl(n,K)$ with $p|n$ should be
omitted (or replaced by the Lie algebras of corresponding general
linear groups).  Then there is always a $G$-equivariant isomorphism
between $\frg$ and $\frg^*$, which transports the Jordan decomposition
in $\frg$ to $\frg^*$.

Earlier work of Kac--Weisfeiler shows that the crucial case to study
is that of a \emph{nilpotent} $\chi \in \frg^*$ (corresponding to some
nilpotent $e \in \frg$).  Here one begins to make connections with the
results on nilpotent orbits summarized above and with related
conjectures arising in Lusztig's work \cite{lu97b,lu98,lu99a,lu99b}.
From now on we consider only the nilpotent case, subject to the above
restrictions on $p$ and $\Phi$.

\subsection{} 

The blocks of $\uxg$ have been determined by Brown and Gordon
\cite{bg01}.  As summarized by Jantzen \cite[C.5]{ja04a}, there is a
natural bijection between the blocks and the ``central characters'',
which in turn are parametrized by the $W$-orbits in $X/pX$ under the
dot action $w \cdot \lam := w(\lam+\rho)-\rho$.  This is a Lie algebra
version of the Linkage Principle.

If $e$ is the nilpotent element corresponding to $\chi$, the component
group $A(e)$ permutes the simple modules in a block.  This action is
understood only in some special cases.

In general the simple modules in a given block are not easy to
parametrize by weights, though each can be obtained as a quotient of
one or more ``baby Verma modules'': these are induced from
one-dimensional modules for a Borel subalgebra $\frb$ satisfying
$\chi(\frb) =0$.  The choice of $\frb$ affects this construction when
$\chi \neq 0$ if $\calb_e$ has more than one irreducible component.

\subsection{}

To make contact with the geometry of $\caln$, we look only at
\emph{regular} blocks: those for which the weight parameters attached
to simple modules lie inside alcoves.  This requires $p>h$ (where $h$
is the Coxeter number).  Jantzen's translation functors then furnish
information about other blocks.

For a regular block of $\uxg$, the work of Bezrukavnikov, Mirkovi\'c,
and Rumynin provides a geometric interpretation.  Under the assumption
that $p>h$, they prove that the number of nonisomorphic simple
modules in the block is equal to the Euler characteristic of
the Springer fiber $\calb_e$: see \cite[5.4.3, 7.1.1]{bmr}.

\subsection{}

The best understood case involves a nilpotent orbit in $\frg^*$
containing some $\chi$ in \emph{standard Levi form}, which means that
the corresponding nilpotent element $e$ is regular in some Levi
subalgebra of a parabolic subalgebra $\frp_I$ of $\frg$ (determined by
a set $I$ of simple roots).  All nilpotent orbits satisfy this
condition for $\frg =\frsl(n,K)$, but in general things get more
complicated.  (See \cite[\S10]{ja98}, \cite[\S2]{ja00},
\cite[D.1]{ja04a}.)

Jantzen has studied simple $\uxg$-modules (and their projective
covers) in considerable detail when $\chi$ has standard Levi form.  In
particular, each simple module can be labelled as $L_\chi(\lam)$ for
one or more $\lam \in X$.  Here $L_\chi(\lam) \cong L_\chi(\mu)$ if
and only if $\mu \in W_I \cdot \lam+pX$, where $W_I$ is the subgroup
of $W$ generated by simple reflections for $\al \in I$ and $w \cdot
\lam :=w(\lam+\rho)-\rho$.

This can be pictured in terms of the alcove geometry of $W_a$, with
the origin of the affine space $E$ taken to be $-\rho$ and the
translations all multiplied by $p$.  Jantzen calls the group $W_p$ in
this setting.  Fixing a weight $\lam$ inside the lowest dominant
alcove, the orbit $W_p \cdot \lam$ under the natural dot action
contains (with periodic repetitions) all weights needed to parametrize
the simple modules in a single regular block.  In fact, it suffices to
work with the $|W|$ alcoves surrounding a single special point such as
$-\rho$.  Then the induced action of $W_I$ on these alcoves identifies
those which correspond to the same simple module.

\subsection{} 

In \cite{ja99b,ja04a}, Jantzen has also studied in depth the case of a
\emph{subregular} $\chi$: its $G$-orbit has dimension $2N-2$, where
$N=|\Phi^+|$.  Only in types $A_n$ and $B_n$ does such an orbit have a
representative in standard Levi form.  But the simple modules in a
regular block of $\uxg$ can be correlated closely with the irreducible
components of $\calb_e$ (here a Dynkin curve), which helps to bypass
the problem of labelling by weights.

\section{Simple modules and left cells}

\subsection{} 

Here we suggest closer connections between the representation theory
discussed in \S3 and the cells in $W_p$.  While our ideas are
speculative, they have some support from computations in special cases
(including unpublished work of Jantzen as well as \cite{hu02}).

Fix a regular block of $\uxg$, with $\chi$ nilpotent, and denote by
$\cals$ a complete set of nonisomorphic simple modules in this block.
As suggested by Bezrukavnikov, this is a candidate for the finite set
$Y$ in Lusztig's formulation discussed in \S2.  If $\chi$ corresponds
to $e \in \frg$, denote by $\call$ the collection of left cells of the
two-sided cell $\Omega$ corresponding in Lusztig's bijection to the
orbit of $e$.  In case the component group $A(e)$ is trivial, the
cardinalities of $\cals$ and $\call$ are expected to be the same:
compare the theorem of \cite{bmr} cited in \S3 with the conjecture
(LC) in \S2.

\begin{conj}
Fix notation as above.
\begin{itemize}
\item[(a)] There is a natural map $\varphi$ from $\cals$ onto $\call$,
whose fibers are the orbits of $A(e)$ in $\cals$.
\item[(b)] A simple module fixed by $A(e)$ maps under $\varphi$ to the
canonical left cell $\Gamma$ in $\Omega$.  $($We call this module
``canonical''.$)$
\end{itemize}
\end{conj}

\subsection{}

The meaning of ``natural'' in part (a) of the conjecture has to be
clarified.  What we have in mind is a simple recipe for assigning
modules to left cells, but it has only been made rigorous in special
cases.  Consider for example the case when $\chi$ has standard Levi
form, so the modules in $\cals$ can be parametrized by weights in a
$W_p$-orbit which lie in alcoves surrounding any given special point
$v \in E$.  Suppose $v$ can be chosen inside the dominant Weyl chamber
in such a way that weights in those surrounding alcoves which lie in
the canonical left cell $\Gamma$ suffice to parametrize $\cals$.  If
$w \in W_p \cong W_a$ labels one of these alcoves, assign the
corresponding simple module to the left cell in $\Omega$ containing
the alcove labelled by $w^{-1}$.  (It would still have to be shown
that this assignment is independent of the choice of the special
point.)

In particular, when $w$ labels the distinguished involution in
$\Gamma$, then $w=w^{-1}$; so the simple module in $\cals$
corresponding to this alcove is assigned to $\Gamma$.  That this
``canonical'' simple module should be fixed by $A(e)$ is suggested by
the parallel discussion in \cite[10.7]{lu89}.

In rank $2$ cases all of this can be observed directly.  But in
general there are serious combinatorial difficulties in working with
the geometry of the cells even in the good case when $\chi$ has
standard Levi form.  The first problem is to locate a suitable special
point $v$.  One might look at the alcove containing the distinguished
involution in the canonical left cell $\Gamma$: this will be the
lowest alcove in $\Gamma$ attached to some special point $v$.  Do the
surrounding alcoves which lie in $\Gamma$ suffice to account for all
simple modules in $\cals$?  In rank $3$, where Gunnells has
constructed pictures of the cells, the evidence about the number of
available alcoves is encouraging.  (But there is one nilpotent orbit
of type $C_3$ which seems to require an alternate choice of special
point.  This orbit has an element in standard Levi form, while the
component group $A(e)$ has order 2.)

\subsection{} 

The highest two-sided cell corresponds to the regular nilpotent
orbit.  Here the related representation theory is quite transparent,
since a regular block has only one simple module (of dimension
$p^N$).

At the other extreme, one can say quite a bit about the lowest
two-sided cell $\Omega$, which corresponds to the zero orbit.  Here
the canonical left cell is just a shifted version of the dominant
chamber, whose geometry is transparent.  The associated representation
theory comes from the group $G$, with simple modules parametrized in
the usual way by highest weights.

Using suggestions of Shi, we can argue as follows.  Start with a
special point for $W_p$ lying in $Q$ such as $v=2(p-1)\rho$; the
surrounding $|W|$ alcoves lie inside the canonical left cell $\Gamma$.
If we write $v = x \cdot (-\rho)$ (with $x$ a translation from $pQ$),
these alcoves are obtained by applying $x$ to the alcoves around
$-\rho$ labelled by the elements $w \in W$, and thus are labelled by
elements $xw$.  Now $x^{-1} \cdot \rho$ lies inside the antidominant
chamber, which is a single left cell of $\Omega$.  Since $W$ acts
simply transitively on the Weyl chambers, we see that the alcoves
labelled by the various $(xw)^{-1}= w^{-1} x^{-1}$ all lie in distinct
Weyl chambers and thus in distinct left cells of $\Omega$.

It is easy to see that the resulting bijection between $\cals$ and
$\call$ is independent of the choice of $v$, since the role of $W$ is
independent of translations by elements of $pQ$.

\subsection{} 

Jantzen's study of the \emph{subregular} case makes it possible to say
something, even though $\chi$ can be chosen to have standard Levi form
only for root systems of type $A_n$ and $B_n$.  In a regular block
there is always an isolated simple module, denoted $L_0$ in
\cite[D.6]{ja99b} and associated with the longest element $w_0$ of
$W$.  This module is characterized in terms of its
``$\kappa$-invariant'' and has a projective cover of smallest possible
dimension.

The dominant alcove $A$ obtained by reflecting the lowest alcove
across its upper wall $\calh$ contains the distinguished involution in
the canonical left cell $\Gamma$; it is the lowest alcove in $\Gamma$
among those sharing the vertex obtained by reflecting $-\rho$ in the
hyperplane $\calh$.  In our framework it is natural to assign the
simple module $L_0$ to $A$ and thus to the left cell $\Gamma$.  (This
is motivated in part by the approach to computing dimensions in
\cite{hu02}, where $\calh$ plays a key role.)  Low rank evidence
indicates that the translate of $A$ attached to the special point
$-\rho$ is in the same $W_I$-orbit as the alcove labelled by $w_0$ in
types $A_n$ and $B_n$.  Here $I$ is the set of simple roots involved
in Jantzen's choice of subregular nilpotent element.

For type $G_2$, there are five simple modules in a regular block,
three of equal dimension being permuted by $A(e) \cong S_3$.  Here the
two-sided cell is finite, with three left cells: the canonical left
cell (to which $L_0$ should be assigned) has 8 elements, while the
others have respectively 8 and 7.  Comparison with Lusztig's model, as
developed by Xi \cite[11.2]{xi94b}, shows that the triple of simple
modules should be assigned to the cell with 7 elements: here the
isotropy group in $S_3$ has order 2.  However, it is unclear for root
systems other than $A_n, B_n, G_2$ how to assign the simple modules
other than $L_0$ to left cells.

\subsection{} 

For a fixed nilpotent orbit, our broader hope is to model Lusztig's
conjectural set-up in full detail.  Besides taking for the finite set
$Y$ the set $\cals$ above, one needs to bring in the action of
$C_G(e)$.  Still missing is a construction (presumably based on
$\calb_e$) of suitable modules which carry compatible actions of
$\frg$ and $F$.

But there is a reasonable prototype in the case $\chi=0$.  Here one
starts with \emph{Weyl modules} $V(\lam)$ with $\lam \in X^+$.  Their
duals are realized as spaces of global sections of line bundles on
$\calb$ (the Springer fiber in this case).  With these modules one has
a Kazhdan--Lusztig theory, conjectured by Lusztig (for $p$ not too
small) to determine simple modules $L(\lam)$ via an alternating sum
formalism with coefficients depending on Kazhdan--Lusztig polynomials
for $W_p$.  In turn $L(\lam)$ factors (by Steinberg's theorem) into a
tensor product of a simple $\uxg$-module and the Frobenius twist of a
simple module for $G$ (which looks like the characteristic $0$ version
if $\lam$ is suitably bounded relative to $p$).

One would like to find a similar construction for all $\chi$.  A
geometric construction of $\frg$-modules using the Springer fiber has
been proposed by Mirkovi\'c--Rumynin \cite{mr01}, but without the
additional features indicated above.

\subsection{}

When $\chi$ is fixed, motivation for correlating simple modules with
left cells comes indirectly from the experimental calculations
reported in \cite{hu02}.  These are reinforced by Jantzen's
unpublished calculations in higher rank cases.  The idea here is that
the geometry of lower boundaries of canonical left cells, together
with the placement of weights in alcoves, should play a key role in
predicting the dimensions (and formal characters) of simple modules.
The experimental evidence also reinforces the suggestion above about
the existence of a tensor product decomposition of Steinberg type.

\subsection{} 

Lusztig's conjectural framework works with a fixed nilpotent orbit or
two-sided cell.  But there is additional motivation for assigning
simple modules to left cells when we compare one orbit with an orbit
in its closure.  When $\psi$ is in the closure of the $G$-orbit of
$\chi$, one expects that a simple $\uxg$-module will ``deform'' into a
not necessarily simple $U_\psi(\frg)$-module.  

On the level of Grothendieck groups, this would imply a recipe for
writing the dimension of the given simple $\uxg$-module as a sum of
dimensions of simple $U_\psi(\frg)$-modules.  In all known cases these
dimension formulas are given uniformly by polynomials in $p$ and the
weight coordinates (compare \cite[\S6]{bmr}).  Experimentation in low
ranks by Jantzen and the author suggests that such decompositions may
be possible in a unique way.

Ostrik proposes that deformation should be studied in the context of
projective covers of simple modules.  He suggests an interpretation of
the process in terms of comparison of Lusztig's equivariant $K$-theory
bases for the two Springer fibers: these bases may be comparable even
when the Springer fibers themselves are not.  Using this viewpoint he
recovers for example our dimension comparisons in the case of type
$G_2$.

In low ranks, the cell pictures related to our hypothetical assignment
of simple modules to left cells show an intriguing correlation with
the computed degeneration in dimension formulas.  But all of this
remains to be placed in a rigorous theoretical setting, beginning
with the process of deformation.  

\bibliographystyle{amsplain}

\begin{thebibliography}{RAGS}

\bibitem{ac04} P.N. Achar, \tit{On the equivariant $K$-theory of the
nilpotent cone in the general linear group}, Represent. Theory
\tbf{8} (2004), 180--211.

\bibitem{as02} P.N. Achar and E.N. Sommers, \tit{Local systems of
nilpotent orbits and weighted Dynkin diagrams}, Represent.
Theory \tbf{6} (2002), 190--201.

\bibitem{be04a} R. Bezrukavnikov, \tit{On tensor categories attached
to cells in affine Weyl groups}, Representation Theory of Algebraic
Groups and Quantum Groups, 69--90, Adv. Stud. Pure Math., \tbf{40},
Math. Soc. Japan, Tokyo, 2004.

\bibitem{be03} \bysame, \tit{Quasi-exceptional sets and
equivariant coherent sheaves on the nilpotent cone},
Represent. Theory \tbf{7} (2003), 1--18.

\bibitem{be02} \bysame, \tit{Perverse sheaves on affine
flags and nilpotent cone of the Langlands dual group},
arXiv:math.RT/0201256.

\bibitem{bmr} R. Bezrukavnikov, I. Mirkovi\'c, D. Rumynin,
\tit{Localization of modules for a semisimple Lie algebra in prime
characteristic}, arXiv:math.RT/0205144, to appear in Ann. of Math.

\bibitem{bo04} R. Bezrukavnikov and V. Ostrik, \tit{On tensor 
categories attached to cells in affine Weyl groups II}, Representation
Theory of Algebraic Groups and Quantum Groups, 101--119,
Adv. Stud. Pure Math., \tbf{40}, Math. Soc. Japan, Tokyo, 2004.

\bibitem{bg01} K.A. Brown and I. Gordon, \tit{The ramification of
centres: Lie algebras in positive characteristic and quantised
enveloping algebras}, Math. Z. \tbf{238} (2001), 733--779.

\bibitem{co02} T. Chmutova and V. Ostrik, \tit{Calculating 
canonical distinguished involutions in the affine Weyl groups},
Experiment. Math. \tbf{11} (2002), 99--117.

\bibitem{go01} I. Gordon, \tit{Representations of semisimple Lie
algebras in positive characteristic and quantum groups at roots of
unity}, pp.~149--167, \tit{Quantum Groups and Lie Theory},
ed. A. Pressley, Proc. Durham 1999, London Math. Soc. Lecture
Note Ser., \tbf{290}, Cambridge Univ. Press, Cambridge, 2001.

\bibitem{hu98} J.E. Humphreys, \tit{Modular representations of simple 
Lie algebras}, Bull. Amer. Math. Soc. (N.S.) \tbf{35} (1998), 105--122.

\bibitem{hu02} \bysame, \tit{Analogues of Weyl's formula for reduced
enveloping algebras}, Experiment. Math. \tbf{11} (2002),
567--573.

\bibitem{ja98} J.C. Jantzen, \tit{Representations of Lie algebras in prime
characteristic}, Notes by Iain Gordon, pp.~185--235, NATO
Adv. Sci. Inst. Ser. C Math. Phys. Sci., 514, \tit{Representation
theories and algebraic geometry (Montreal, 1997)}, Kluwer Acad. Publ.,
Dordrecht, 1998.

\bibitem{ja99b} \bysame, \tit{Subregular nilpotent representations of Lie 
algebras in prime characteristic}, Represent. Theory \tbf{3}
(1999), 153--222.

\bibitem{ja00} \bysame, \tit{Modular representations of reductive 
Lie algebras}, J. Pure Appl. Algebra \tbf{152} (2000),
133--185.

\bibitem{ja04a} \bysame, \tit{Representations of Lie algebras in
positive characteristic}, Representation Theory of Algebraic Groups
and Quantum Groups, 175--218, Adv. Stud. Pure Math., \tbf{40},
Math. Soc. Japan, Tokyo, 2004.

\bibitem{ja04b} \bysame, \tit{Nilpotent orbits in representation theory},
pp.~1--211, \tit{Lie Theory}, ed.~J.-P. Anker and B. Orsted, Progr. 
Math., vol.~228, Birkh\"auser, Boston, 2004.

\bibitem{lu83} G. Lusztig, \tit{Some examples of square integrable
representations of semisimple $p$-adic groups}, Trans. Amer. Math. Soc.
\tbf{277} (1983), 623--653.

\bibitem{lu85} \bysame, \tit{Cells in affine Weyl groups}, 
Algebraic Groups and Related Topics (Kyoto/Nagoya, 1983), 255--287,
Adv. Stud. Pure Math., \tbf{6}, Math. Soc. Japan, Tokyo, 1985.

\bibitem{lu87a} \bysame, \tit{Cells in affine Weyl groups II}, 
J. Algebra \tbf{109} (1987), 536--548.

\bibitem{lu87b} \bysame, \tit{Cells in affine Weyl groups III}, 
J. Fac. Sci. Univ. Tokyo Sect IA Math. \tbf{34} (1987), 223--243.

\bibitem{lu89} \bysame, \tit{Cells in affine Weyl groups IV}, 
J. Fac. Sci. Univ. Tokyo Sect IA Math. \tbf{36} (1989),
297--328.

\bibitem{lu97a} \bysame, \tit{Cells in affine Weyl groups and tensor
categories}, Adv. Math. \tbf{129} (1997), 85--98.

\bibitem{lu97b} \bysame, \tit{Periodic $W$-graphs}, Represent. Theory
  \tbf{1} (1997), 207--279.

\bibitem{lu98} \bysame, \tit{Bases in equivariant $K$-theory},
Represent. Theory \tbf{2} (1998), 298--369.

\bibitem{lu99a} \bysame, \tit{Subregular nilpotent elements and bases 
in $K$-theory}, Canad. J. Math. \tbf{51} (1999), 1194--1225.

\bibitem{lu99b} \bysame, \tit{Bases in equivariant $K$-theory, II},
Represent. Theory \tbf{3} (1999), 281--353.

\bibitem{lx88} G. Lusztig and N. Xi, \tit{Canonical left cells in affine
Weyl groups}, Adv. Math. \tbf{72} (1988), 284--288.

\bibitem{mr01} I. Mirkovi\'c and D. Rumynin, \tit{Geometric representation 
theory of restricted Lie algebras}, Transform.  Groups
\tbf{6} (2001), 175--191.

\bibitem{os00} V. Ostrik, \tit{On the equivariant $K$-theory of the
nilpotent cone}, Represent. Theory \tbf{4} (2000), 296--305.

\bibitem{sh86} Jian Yi Shi, \tit{The Kazhdan--Lusztig cells in 
certain affine Weyl groups}, Lect. Notes in Math. 1179,
Springer-Verlag, Berlin, 1986.

\bibitem{sh88} \bysame, \tit{A two-sided cell in an affine Weyl
group, II}, J. London Math. Soc. (2) \tbf{37} (1988), 253--264.

\bibitem{so01} E. Sommers, \tit{Lusztig's canonical quotient and generalized
duality}, J. Algebra \tbf{243} (2001), 790--812.

\bibitem{xi94a} N. Xi, \tit{The based ring of the lowest two-sided cell
of an affine Weyl group. II}, Ann. Sci. \'Ecole Norm. Sup. (4) 
\tbf{27} (1994), 47--61.

\bibitem{xi94b} \bysame, \tit{Representations of affine Hecke algebras}, 
Lect. Notes in Math. 1587, Springer-Verlag, Berlin, 1994.

\bibitem{xi02} \bysame, \tit{The based ring of two-sided cells of affine
Weyl groups of type $\widetilde{A}_{n-1}$}, Mem. Amer. Math. Soc.
\tbf{157} (2002), no.~749.

\end{thebibliography}

\makeatletter \renewcommand{\@biblabel}[1]{\hfill#1.}\makeatother

\end{document}